\documentclass[12pt]{article}
\usepackage{amsmath}

\usepackage{amssymb}

\newsymbol \blackbox 1004
\newcommand{\eh}{\hfill}\newlength{\sperr}

\newenvironment{proof}{{\settowidth{\sperr}{\bf\rm Proof}%
\par\addvspace{0.3cm}\noindent\parbox[t]{1.3\sperr}
{\bf\rm P\eh r\eh o\eh o\eh f\eh }%
}}{\nopagebreak\mbox{} $\blackbox$\par\addvspace{0.3cm}}

\def\a{\alpha}

\def\g{\gamma}

\def\s{\sigma}
\def\t{\theta}

\def\ov{\overline}
\def\vk{\varkappa}
\def\wt{\widetilde}
\def\BC{{\mathbb C}} % \mathbb
\def\BR{{\mathbb R}} % \mathbb

\def\im{{\rm Im\ }}

\newtheorem{Pa}{Paper}[section]
\newtheorem{Tm}[Pa]{{\bf Theorem}}
\newtheorem{La}[Pa]{{\bf Lemma}}
\newtheorem{Cy}[Pa]{{\bf Corollary}}

\newtheorem{Pn}[Pa]{{\bf Proposition}}

\date{}

\title{Completion problems and scattering problems for Dirac type
differential equations
with singularities}

\author{Bernd Fritzsche, Bernd Kirstein, and Alexander Sakhnovich }

\begin{document}
\maketitle

\begin{abstract}
A completion problem to recover a rational matrix function which is $j$-unitary on the line  is treated.
A Dirac type system with singularities on the semiaxis is recovered explicitly by its left
reflection coefficient. The close connection between these two problems is discussed.
\end{abstract}

% main text
\section{Introduction} \label{intro}
\setcounter{equation}{0}
Completion problems is an interesting domain closely connected
with  moment and interpolation problems (see, for instance, \cite{Ar,Ar6,AFK0,AFK,B1,DFK,FK1,FK2}).
In this paper we shall show
that completion problems underlie
 explicit solutions of some classical and non-classical inverse scattering problems   too.
We  call an $m \times m$ rational matrix function $W(\lambda)$,
that tends to $I_m$ at infinity, {\it weakly $j$-elementary} if it
is $j$-unitary on the real axis:
 \\
a) $\, W(\lambda)^*jW(\lambda)=j \,$ for $\lambda$ on the real
line $\BR$.
\\ Here $I_m$ is the $m \times m$ identity matrix,$\, W^*$ stands
for the conjugate transpose of matrix $W$, and $\, j=j^*=j^{-1}$.
We call $W$ {\it regular} in $\BC_+$ if: \\ b) $\, W(\lambda)$ has
no singularities in the open upper halfplane $\BC_+$. \\ Notice
that a matrix function $W(\lambda)$, that tends to $I_m$ at
infinity, is called $j$-{\it elementary rational factor in} $
\BC_+$ if: \\ c) $\, W(\lambda)^*jW(\lambda) \geq j \,$ for
$\lambda \in \BC_+$,  $\, W(\lambda)^*jW(\lambda) \leq j \,$ for
$\lambda \in \BC_-$,
\\ but the weak requirement a) deals with the axis only. Partition
now $W$ into four blocks $W= \{ W_{kj} \}_{k,j=1}^2$, where
$W_{11}$ is an $m_1 \times m_1$ matrix.

Our main completion problem is to recover a weakly $j$-elementary and regular function $W$ by
\begin{equation} \label{0.3}
R(\lambda)= W_{21}(\lambda)W_{11}(\lambda)^{-1}.
\end{equation}

The $j$-elementary and weakly $j$-elementary rational factors in  $\BC_+$ prove to be closely connected
\cite{GKS1,SaA1,SaAp} with the
fundamental solution of the self-adjoint Dirac type (also called ZS or AKNS)
system, which is a classical matrix differential equation:
\begin{equation}       \label{1.1}
\frac{d}{dx}u(x, \lambda  )=i(\lambda  j+jV(x))u(x, \lambda  )
\quad (x \geq 0),
\end{equation}
where $m=m_1+m_2$,
\begin{equation}   \label{1.2}
j = \left[
\begin{array}{cc}
I_{m_1} & 0 \\ 0 & -I_{m_2}
\end{array}
\right], \hspace{1em} V= \left[\begin{array}{cc}
0&v\\v^{*}&0\end{array}\right].
 \end{equation}
Dirac type systems are very well-known in mathematics and applications (see, for instance,
\cite{CS,CG,DI,K,LS,SaA3,SaL3,YL} and references therein). The name
ZS-AKNS is connected with the fact that system (\ref{1.1}) is an
auxiliary linear system for many important nonlinear integrable
wave equations and as such it was studied, for instance, in
\cite{AKNS,AS,FT,GSS,ZS}.

For simplicity we shall further consider $j$ of the form (\ref{1.2}).
First we solve in this paper the above mentioned  completion problem. Then we treat
Dirac type systems on the semiaxis with left reflection coefficient $R_L$ satisfying (\ref{0.3}).
We consider $m_1 \times m_2$ rectangular and slowly decaying potentials $v$ with singularities that are essential
in the study of the multicomponent nonlinear equations and soliton-positon interactions \cite{BKo,BM,Kur,M1,Sch}.
Explicit solutions of the direct and inverse problems that generalize and develop earlier results
from \cite{AGKS,AD,GKS6,SaAp} are obtained.
Various notions from
the mathematical system theory will be widely used, and therefore an Appendix with  basic notions and results
from system theory
 is added to the paper.

It is well-known \cite{SaL1,AG0} that any  rational matrix function which is $j$-unitary on $\BR$ and
tending to $I_m$ at infinity
can be presented in the form of the transfer matrix function
\begin{equation} \label{0.1'}
W( \lambda)=I_m+ij \Lambda_0^*S_0^{-1}(\lambda I_n- \alpha)^{-1}
\Lambda_0,
\end{equation}
where $\alpha$ and $S_0$ are square $n \times n$ matrices,
$\Lambda_0$ is $n \times m$ matrix, $n$ is the McMillan degree of
$W$, $S_0=S_0^*$, and the identity
\begin{equation} \label{0.2'}
\alpha S_0-S_0 \alpha^*=i \Lambda_0 j \Lambda_0^*
\end{equation}
holds. The transfer matrix functions of the special form
(\ref{0.1'}), (\ref{0.2'}) have been introduced and studied in
\cite{SaL1}-\cite{SaL3} (see also the references therein).
Regularity of $W$ means that the spectrum of $\a$ in the
realization  (\ref{0.1'}) belongs to the closed lower halfplane
$\ov \BC_-$. When $S_0>0$ we can substitute $\a$, $S_0$,  and
$\Lambda_0$ in (\ref{0.1'}), (\ref{0.2'}) by $ S_0^{- \frac{1}{2}}
\a S_0^{ \frac{1}{2}}$, $I_n$,  and $S_0^{- \frac{1}{2}}
\Lambda_0$, respectively, i.e., without loss of generality we can
assume that $S_0=I_n$ in this case. Relations (\ref{0.1'}),
(\ref{0.2'}) with $S_0=I_n$ take the form
\begin{equation} \label{0.1}
W( \lambda)=I_m+ij \Lambda_0^*(\lambda I_n- \a)^{-1} \Lambda_0,
\quad \a - \a^*=i \Lambda_0 j \Lambda_0^* .
\end{equation}
Notice that
the matrix function in (\ref{0.1}) is a well-known Liv\v{s}ic-Brodskii characteristic
function \cite{L}. Important applications of its multidimensional generalization
to the spectral and scattering theory one can find in \cite{B2,V} and the references therein.

\begin{Tm} \label{Tm0.1}
$\quad (i)$ The weakly $j$-elementary, regular  in  $ \BC_+$ rational matrix function $W$
tending to $I_m$ at infinity
is uniquely recovered by the $m_2 \times m_1$ matrix function $R$ satisfying (\ref{0.3}).
This matrix function $R$ is contractive on $\BR$ and tends to zero at infinity. \\
$\quad (ii)$ Moreover, by each rational, contractive on $\BR$ and tending to zero at infinity $m_2 \times m_1$
matrix function $R$ we can recover a unique weakly $j$-elementary and regular  in  $ \BC_+$ rational
matrix function $W$
tending to $I_m$ at infinity,
of the same McMillan degree as $R$, and
such that (\ref{0.3}) holds.
The matrix function  $W$ is recovered in the following way. Consider a minimal
realization
\begin{equation}\label{2.15}
R(\lambda )= C(\lambda I_n- A)^{-1}B,
\end{equation}
where $n$ is the McMillan degree of $R$, i.e., the order of $A$,
and choose the $n \times n$ Hermitian solution $X=X^*$ of the Riccati equation
\begin{equation}\label{Ric}
i(XA-A^*X)=C^*C+XBB^*X
\end{equation}
such that
\begin{equation}\label{2.16}
\s(A+iBB^*X) \subset \ov \BC_- .
\end{equation}
Here $\s$ means spectrum and   ${\ov \BC}_-$ is the closed lower halfplane.  Next put
\begin{equation}\label{0.4}
\a=A+i BB^*X, \quad S_0=X^{-1}, \quad \g_1=
B, \quad  \g= -i S_0 C^*,
\end{equation}
\[ \Lambda_0=[\g_1, \quad  \g]. \]
Finally substitute (\ref{0.4}) into (\ref{0.1'}) to construct $W$. \\
$\quad (iii)$ If $W$ is a $j$-{\it elementary rational regular factor in} $ \BC_+$, then
$R$ is contractive in $ \BC_+$, and vice versa if the conditions of $(ii)$ hold and
$R$ is contractive in $ \BC_+$,
then $X$ in (\ref{Ric})-(\ref{0.4}) is strictly positive and the recovered $W$
is a $j$-{\it elementary rational regular factor in} $ \BC_+$.
\end{Tm}
The completion problem considered in $(iii)$ has been in a different way solved in
\cite{Ar,DeD} (see also \cite{GKS6} for the subcase of (iii), where $m_1=m_2$).
When $\s(\a) \subset \BR$ the matrix function $W$ belongs to the class
of {\it Arov-singular} $j$-inner matrix functions \cite{Ar6,AFK0,AFK}. The case $\s(\a) \subset \BR$ corresponds to
system (\ref{1.1}) with the slow-decaying $n$-positon type potentials $v$ \cite{SaAp}.

To prove Theorem \ref{Tm0.1} we shall need the following generalization (for the case
of the sign-indefinite solutions $X$ of the Riccati equations) of Theorem 4.2 given in the interesting
paper \cite{AD} which is related to our work.
\begin{La} \label{La0.2} Suppose the $m_2 \times m_1$
matrix function $R$ is rational, contractive on $\BR$ and tending to zero at infinity,
and (\ref{2.15}) is its minimal realization. Then there is a unique Hermitian solution $X$
of the Riccati equation  (\ref{Ric}) such that (\ref{2.16}) holds. This solution $X$ is always invertible.
It is also strictly positive if and only if $R$ is  contractive in $\BC_+$.
\end{La}
This lemma is closely connected with several results from \cite{LR}.

Lemma \ref{La0.2} and Theorem \ref{Tm0.1} will be proved in the  Section 2.
Explicit solutions of the direct scattering problem are given in Section 3
whereas explicit solutions of the inverse scattering problem are given in Section 4.
An appendix contains a collection of  necessary notions and results from  system theory.
%%%%%%%%%%%%%%%%%%%%%%%%%%%%%%%%%%%%%%%%%%%%%%%%%%%%%%%%%%%%%%%%%%%%%%%%%%%%%%%%%%%%%%%%%%%%%%%%%%%%%%%%
%%%%%%%%%%%%%%%%%%%%%%%%%%%%%%%%%%%%%%%%%%%%%%%%%%%%%%%%%%%%%%%%%%%%%%%%%%%%%%%%%%%%%%%%%%%%%%%%%%%%%%%%
\section{Completion problem} \label{compl}
\setcounter{equation}{0}
\begin{proof} of Lemma \ref{La0.2}.
Denote by $\ov \lambda$ complex value conjugated to $\lambda$ and
by $\mathcal M^*$, for some set $\mathcal M$ of complex values,
the set of complex values conjugated to those in $\mathcal M$. By
Theorem 21.1.3 in \cite{LR} there exists a Hermitian solution $X$
of the Riccati equation (\ref{Ric}). (Indeed, the matrix function
$R(\ov \lambda)^*=B^*(\lambda I_n-A^*)^{-1}C^*$ is  contractive on
$\BR$ simultaneously with $R(\lambda)$ and equation (21.1.6) in
\cite{LR} written down for $R(\ov \lambda)^*$ coincides with
(\ref{Ric}).) We introduce the characteristic state space matrix
\cite{GR}
\begin{equation}\label{L1}
K=\left[
\begin{array}{cc}
A & BB^* \\ C^*C & A^*
\end{array}
\right].
\end{equation}
If $X$
satisfies (\ref{Ric})  then we have
\begin{equation}\label{L2}
T^{-1}KT=\left[
\begin{array}{cc}
A+i BB^* X & BB^* \\ 0 & (A+i BB^* X)^*
\end{array}
\right], \quad T=
\left[
\begin{array}{cc}
I_n & 0 \\ i X & I_n
\end{array}
\right].
\end{equation}
By (\ref{L2}) we obtain
\begin{equation}\label{L3}
\s(K)\setminus \BR= \mathcal M \cup  \mathcal M^*,
\end{equation}
where  $\mathcal M$ is the subset of the eigenvalues of $K$ in
$\BC_-$. The subset $\mathcal M$  has the evident property: if
$\mu \in \mathcal M$, then $\ov \mu \not\in \mathcal M$. Moreover,
according to (\ref{L3}) $\mathcal M$ is a maximal set of non-real
eigenvalues of $K$ subject to this property. Now we can apply
Theorem 7.4.2 in \cite{LR} to derive the existence and the
uniqueness of the Hermitian solution $X$ of the Riccati equation
(\ref{Ric}) such that $\s(A+iBB^*X) \setminus \BR= \mathcal M$
from the maximality of $\mathcal M$ and the existence of some
Hermitian solution of (\ref{Ric}). (Though normalizations in
Theorem 7.4.2 in \cite{LR} are chosen for the right halfplane they
are easily changed to suit $\BC_-$.) Use again (\ref{L2}) to see
that (\ref{2.16}) yields $\s(A+iBB^*X) \setminus \BR= \mathcal M$.
Therefore a Hermitian matrix $X$ satisfying (\ref{Ric}) and
(\ref{2.16}) is unique too. Finally notice that by Theorem 21.2.1
in \cite{LR} any Hermitian $X$ satisfying (\ref{Ric}) is
invertible, and this $X$ is strictly positive if and only if
$R(\ov \lambda)^*$ has no poles in the lower halfplane. In other
words $X$ is strictly positive if and only if $R( \lambda)$ has no
poles in $\BC_+$, i.e., if and only if $R( \lambda)$ is
contractive in $\BC_+$.
\end{proof}
\begin{proof} of Theorem \ref{Tm0.1}.
According to (\ref{0.1'}), if $\lambda \to \infty$  then we get
$W_{11}(\lambda  ) \to I_{m_1}$ and $W_{21}(\lambda  ) \to 0$. So
$R$ of the form (\ref{0.3}) tends to zero. As
$W(\lambda)^*jW(\lambda) = j$ on $\BR$ we have
\begin{equation} \label{n1.2}
W_{11}(\lambda  )^*W_{11}(\lambda  )-W_{21}(\lambda
)^*W_{21}(\lambda  ) = I_{m_1},
\end{equation}
in particular. In view of (\ref{n1.2}), $R$ of the form (\ref{0.3}) is contractive on $\BR$.
Thus the necessary requirements for $R$ are proved.

Suppose (\ref{2.15}) is a minimal realization of some
matrix function $R$ contractive  on $\BR$. By Lemma \ref{La0.2}
 there is a Hermitian $X$ satisfying (\ref{Ric}) and (\ref{2.16}).
Therefore the procedure to recover $W$ in theorem is well defined. Moreover, according to
(\ref{0.4}) the Riccati equation (\ref{Ric}) is equivalent to (\ref{0.2'}), i.e., (\ref{0.2'}) is valid.
In view of (\ref{0.1'}) and  (\ref{0.4}) we have
\[
W_{11}(\lambda)=I_{m_1}+i\g_1^*S_0^{-1}(\lambda  I_n - \a)^{-1}
\gamma_1, \quad W_{21}(\lambda)=
\]
\begin{equation} \label{nn1.1}
=-i\g^*S_0^{-1}(\lambda  I_n - \a)^{-1} \gamma_1.
\end{equation}
Using formula (\ref{app4}) in the Appendix, from the first equality in (\ref{nn1.1}) we get
\begin{equation} \label{nn1.2}
W_{11}(\lambda)^{-1}=I_{m_1}-i\g_1^*S_0^{-1}(\lambda  I_n -
\t)^{-1} \gamma_1, \quad \t=\a-i \gamma_1\gamma_1^*S_0^{-1}=A.
\end{equation}
From (\ref{nn1.1}) and (\ref{nn1.2})
we obtain
\[
W_{21}(\lambda)W_{11}(\lambda)^{-1}= -i\g^*S_0^{-1}(\lambda  I_n
-\a)^{-1} \gamma_1
\]
\begin{equation} \label{fn1.3}
-\g^*S_0^{-1}(\lambda  I_n -\a)^{-1} \gamma_1
\g_1^*S_0^{-1}(\lambda  I_n - \t)^{-1} \gamma_1.
\end{equation}
Taking into account the identity
\[
\gamma_1\gamma_1^*S_0^{-1}=i(\t- \a)=-i\big( (\lambda I_n-
\t)-(\lambda I_n- \a) \big),
\]
we can rewrite (\ref{fn1.3}) as
\begin{equation} \label{n1.3}
W_{21}(\lambda)W_{11}(\lambda)^{-1}= - i\g^*S_0^{-1}(\lambda  I_n
- \t)^{-1} \gamma_1 ,
\end{equation}
and the right-hand side of (\ref{n1.3}) coincides with $R$.
In other words, equality (\ref{0.3}) is valid.
Moreover, the McMillan degree of the recovered $W$ equals $n$, i.e., the realization (\ref{0.1'}) is
minimal. Indeed,  the right-hand
side of (\ref{0.3})  has its McMillan degree less or equal to the McMillan degree of $W$,
and so the McMillan degree of $R$ is less or equal to the McMillan degree of $W$. Therefore the degrees of $W$ and
$R$ are equal.
The formulas (\ref{0.1'}) and (\ref{0.2'})
yield the widely used equality
\begin{equation} \label{n1.1}
W(\lambda  )^*jW(\lambda  )=j+i(\ov{\lambda}- \lambda) \Lambda_0^*
(\ov{\lambda} I_n -\alpha^*)^{-1}S_0^{-1}(\lambda I_n
-\alpha)^{-1}\Lambda_0.
\end{equation}
By  (\ref{n1.1}) the function $W$ is a weakly $j$-elementary factor. According to (\ref{2.16}) and (\ref{0.4}),
we have
$\s(\a) \subset \ov \BC_-$, and so $W$ of the form (\ref{0.1'}) has no singularities in $\BC_+$, i.e.,
$W$ is regular.
%%%%%%%%%%%%%%%%%%%%

To prove uniqueness
suppose that there is another weakly  $j$-elementary regular rational factor $\wt W$ in  $ \BC_+$
of the  McMillan degree $n$ and such that the equality $R= \wt W_{21} \wt W_{11}^{-1}$ is valid.
Then $\wt W$ admits a minimal realization such that
\begin{equation} \label{n0.1}
\wt W( \lambda)=I_m+ij \wt \Lambda_0^* \wt S_0^{-1}(\lambda I_n-
\wt \a)^{-1} \wt \Lambda_0,
\end{equation}
\begin{equation} \label{n0.1'}
\wt \a \wt S_0- \wt S_0 \wt \a^*=i \wt \Lambda_0 j \wt
\Lambda_0^*, \quad \wt \Lambda_0=[\wt \g_1, \quad \wt \g], \quad
\s( \wt \a) \subset \ov \BC_-.
\end{equation}
As $R= \wt W_{21} \wt W_{11}^{-1}$ the matrix function $R$ admits another minimal realization, too:
\begin{equation} \label{n1.4}
R(\lambda )= \wt C(\lambda  I_n- \wt A)^{-1}\wt B,
\end{equation}
where
\begin{equation} \label{n1.5}
\wt A= \wt \a-i \wt \g_1 \wt  \g_1^* \wt  S_0^{-1}, \quad \wt B= \wt  \g_1,
\quad  \wt C= -i \wt \g^* \wt S_0^{-1}
\end{equation}
(see (\ref{n1.3})). Therefore
there is a similarity transformation  matrix $s$ such that:
\begin{equation} \label{n1.7}
\wt A= s A s^{-1}, \quad \wt B=s  B, \quad  \wt C=Cs^{-1}.
\end{equation}
Taking into account (\ref{Ric}), (\ref{2.16}), and (\ref{n1.7}), we see that the matrices $\wt X=(s^*)^{-1}Xs^{-1}$
and $\wt A+i \wt B \wt B^* \wt X=s \a s^{-1}$
satisfy the relations
\begin{equation}\label{Ric'}
i(\wt X \wt A-\wt A^*\wt X)=\wt C^*\wt C+\wt X \wt B \wt B^* \wt X,
\quad \s(\wt A+i \wt B \wt B^* \wt X) \subset \ov \BC_- .
\end{equation}
By (\ref{n0.1'}) and (\ref{n1.5}) we see that $\wt X=\wt S_0^{-1}$ satisfies (\ref{Ric'}), too.
But, according to Lemma \ref{La0.2}, the Hermitian solution $\wt X$ satisfying (\ref{Ric'}) is unique, i.e.,
\begin{equation}\label{nn1.4}
\wt S_0=sX^{-1}s^*=sS_0s^*.
\end{equation}
From (\ref{0.4}), (\ref{n1.5}), (\ref{n1.7}), and (\ref{nn1.4}) it follows that
\begin{equation}\label{nn1.5}
\wt \Lambda_0=s \Lambda_0, \quad \wt \a = s \a s^{-1}.
\end{equation}
By (\ref{0.1'}), (\ref{n0.1}), (\ref{nn1.4}), and (\ref{nn1.5}) we derive
$W= \wt W$, and so the uniqueness of $W$ is proved.

Finally, if $W$ is $j$-elementary we get
\[
W_{11}(\lambda  )^*W_{11}(\lambda  )-W_{21}(\lambda
)^*W_{21}(\lambda  ) \geq I_{m_1} \quad (\lambda \in \BC_+),
\]
i.e., $R^*R \leq I_{m_1}- (W_{11}^*)^{-1}W_{11}^{-1}$.  Thus the first statement in (iii) is proved.
On the other hand, if $R$ is contractive in $ \BC_+$, then by Lemma \ref{La0.2} we have $S_0=X^{-1}>0$.
Use now (\ref{n1.1}) to show that $W$ is $j$-elementary. In this way, the second statement in (iii)
is proved, too.
\end{proof}
From (\ref{0.2'}) it follows that
\begin{equation}\label{nn1.6}
\t S_0 - S_0 \t^*=-i \Lambda_0 \Lambda_0^*, \quad \t=\a-i
\gamma_1\gamma_1^*S_0^{-1}, \quad \Lambda_0 = [\g_1 \quad \g ].
\end{equation}
Similar to the formula (\ref{n1.1}) we obtain for the matrix function
\begin{equation}\label{nn1.7}
{\bf S}(\lambda)=I_m - i  \Lambda_0^* S_0^{-1}(\lambda I_n - \t
)^{-1} \Lambda_0
\end{equation}
the equality
\begin{equation} \label{n1.8}
{\bf S}(\lambda  )^* {\bf S}(\lambda  )=I_m+i(\lambda -
\ov{\lambda}) \Lambda_0^* (\ov{\lambda} I_n
-\t^*)^{-1}S_0^{-1}(\lambda I_n -\t)^{-1}\Lambda_0.
\end{equation}
Therefore Theorem \ref{Tm0.1} yields the  solution of the so called "minimal unitary completion problem".
\begin{Cy} \label{Cy1.1} Suppose that the $m_2 \times m_1$
matrix function $R$ is rational, contractive on $\BR$ and tending
to zero at infinity, and that (\ref{2.15}) is its minimal
realization. Then the matrix function ${\bf S}(\lambda  )$ given
by the formula (\ref{nn1.7}) and by the relations
(\ref{Ric})-(\ref{0.4}) is a unitary completion of $R$,  where $R$
is its left lower block, and has the same McMillan degree as $R$.
\end{Cy}
The minimal unitary completion problem for the strictly contractive matrix functions has been solved
in \cite{GR} (see also \cite{LR}). It will be shown in the next section (see formula (\ref{n3.7})) that  one of the two
unitary factors  of the four block scattering matrix function for the Dirac type system takes the form
(\ref{nn1.7}).

%%%%%%%%%%%%%%%%%%%%%%%%%%%%%%%%%%%%%%%%%%%%%%%%%%%%%%%%%%%%%%%%%%%%%%%%%%%%%%%%%%%%%%%%%%%%%%%%
%%%%%%%%%%%%%%%%%%%%%%%%%%%%%%%%%%%%%%%%%%%%%%%%%%%%%%%%%%%%%%%%%%%%%%%%%%%%%%%%%%%%%%%%%%%%%%%%%%

\section{Direct scattering problem for Dirac type systems} \label{direct}
\setcounter{equation}{0}

Classical scattering results for  Dirac type systems, the
history, and literature one can find in \cite{FT}. For the more
recent, matrix, and non-classical scattering results, see
\cite{AKM,APP,BC,BDT,CS,Zh} and the references therein. Explicit solutions
of the direct and inverse scattering problems and nonlinear equations
(soliton, positon, negaton  solutions, in particular) are of special
interest.

We shall consider system (\ref{1.1}) with the so called
pseudo-exponential $m_1 \times m_2$ potential $v$ of the form
\begin{equation} \label{1.3}
v(x)=-2i \gamma _{1}^{\, *}e^{ix \alpha ^{*}} S (x)^{-1}e^{ix
\alpha }\g,
\end{equation}
where
the $n \times n$ matrix function $S$ is given by the formula
\begin{equation}   \label{1.4}
S (x)=S_0+ \int_{0}^{x} \Lambda(t) \Lambda (t)^{*}dt, \hspace{1em}
\Lambda (x)= \left[  e^{-ix \alpha } \gamma_{1}  \qquad
e^{ix \alpha } \gamma
\right],
\end{equation}
$\a$, $\g_1$, and $\g$ are $n \times n$, $n \times m_1$, and $n
\times m_2$ parameter matrices, and the following matrix identity is satisfied:
\begin{equation} \label{1.5}
\a S_0-S_0 \a^*=i(\g_1 \g_1^*- \g \g^*) \qquad (S_0=S_0^*).
\end{equation}
Various direct and inverse problems for system (\ref{1.1}),
(\ref{1.3}) can be solved in terms of the parameter matrices
explicitly. These problems for rapidly decaying
subclasses of $v$ of the
form  (\ref{1.3}) with $\s(\a) \subset \BC_-$ were studied, in particular, in \cite{AG,AGKS,AKM,AD}.
When $\s(\a) \cap \BR \not= \emptyset$ the potential $v$ is slowly decaying
and is not necessarily summable on any semiaxis $(c, \infty)$ \cite{GKS6}.
When $\s(\a) \subset \BR$ the scattering matrix for system (\ref{1.1}),
(\ref{1.3}) on the whole axis equals $I_m$, and in this case we have the supertransparent
($n$-positon) potentials \cite{SaAp}.
The spectral theory of system (\ref{1.1}), (\ref{1.3}) on the
semi-axis $[0, \infty)$, when $m_1=m_2$ and $S_0=I_n$, was studied
in \cite{GKS1,GKS4,GKS6} and other papers of the same authors
cited therein. In particular, it was shown in Theorem 5.4 in
\cite{GKS6} that without loss of generality one can assume:
\begin{equation}\label{1.6}
\s(\a) \subset \overline{\BC}_-, \quad
{\mathrm{span}}\bigcup_{l=0}^{n-1} {\mathrm{Im}}\,\a^l \g_1=
{\mathbb C}^n, \quad {\mathrm{span}}\bigcup_{l=0}^{n-1}
{\mathrm{Im}}\,\a^l \g= {\mathbb C}^n,
\end{equation}
where Im  is image, and the pair $\a$, $\g$ that satisfies the third
equality in (\ref{1.6}) is called full range or controllable. We shall
call the set of parameter matrices $\a$, $S_0$, $\g_1$, and $\g$
{\it admissible} if relations (\ref{1.5}) and (\ref{1.6}) hold.
Analogously to \cite{GKS6} the class of potentials $v$ of the form
(\ref{1.3}) {\it determined} via (\ref{1.4}) by a set of admissible parameter
matrices will be denoted by the acronym PE. In this and the next sections
we generalize a part of  results in \cite{GKS6} to include the cases $m_1 \not= m_2$
and potentials with singularities. (Singularities occur on the semiaxis when
the inequality $S_0>0$ does not hold.)

Similar to the  case $m_1=m_2$ the special $m \times m_1$ and $m \times m_2$
solutions $Y$ and $Z$,
respectively, of system (\ref{1.1}) we determine by
the relations
\begin{equation}\label{1.7}
 Y(x,\lambda )=e^{ix\lambda}\left[ \begin{array}{c} I_{m_1} \\ 0 \end{array}
  \right] +o(1)
\quad (x\to\infty ),\qquad Z(0,\lambda )=\left[ \begin{array}{c} 0
\\ I_{m_2} \end{array} \right] ,
\end{equation}
where $\lambda \in \BR$. In terms of these special solutions the
so called left and right transmission and reflection coefficients
of system (\ref{1.1}) can be expressed via
\begin{eqnarray}
&&T_L(\lambda ):=Y_1(0,\lambda )^{-1},\quad
R_L(\lambda):=Y_2(0,\lambda )\, Y_1(0,\lambda )^{-1},\label{1.8}\\
\noalign{\vskip6pt} && R_R(\lambda ):=
\Bigl(\lim_{x\to\infty}e^{-ix\lambda}Z_1(x,\lambda )\Bigr)
\Bigl(\lim_{x\to\infty}e^{ix\lambda}Z_2(x,\lambda )\Bigr)^{-1},
\label{1.9}\\ \noalign{\vskip6pt}
&&T_R(\lambda):=\Bigl(\lim_{x\to\infty}e^{ix\lambda}Z_2(x,\lambda
)\Bigr)^{-1}.\label{1.10}
\end{eqnarray}
Here the $m_1 \times m_1$ matrix $Y_1$ and the $m_1 \times m_2$ matrix $Y_2$ are upper
and lower blocks of $Y$, respectively. Analogously, the $m_1 \times m_2$ matrix $Z_1$
and the $m_2 \times m_2$ matrix $Z_2$ are upper
and lower blocks of $Z$, respectively.
The functions $T_L$ and $T_R$ are called the {\em left and right transmission
coefficients} and $R_L$ and $R_R$ the {\em left and right reflection
coefficients},  respectively.
Then the four block scattering matrix function ${\mathcal{S}}$ is
given by the equality
\begin{equation} \label{1.11}
{\mathcal S}(\lambda )=    \left[
\begin{array}{cc}
T_L(\lambda ) & R_R(\lambda ) \\ \noalign{\vskip6pt} R_L(\lambda )
& T_R(\lambda )
\end{array}
\right].
\end{equation}
To construct $Y$ and $Z$ for system (\ref{1.1}), (\ref{1.3}) we
shall need some preliminaries. The $m \times m$ solution $u(x,
\lambda)$ of system (\ref{1.1}), (\ref{1.3})  of the form
\begin{equation} \label{1.12}
u(x, \lambda)=w_{\a}(x, \lambda)e^{ix \lambda j},
\end{equation}
where
\begin{equation} \label{1.13}
w_{\a}(x, \lambda)=I_m+ij \Lambda(x)^*S(x)^{-1}(\lambda I_n-
\a)^{-1} \Lambda(x),
\end{equation}
and where $n$ is the order of the parameter matrix $\a$ (and the dimension
of the so called state space), was constructed in Theorem 3 in
\cite{SaA1}. See also \cite{GKS1}, where the case $m_1=m_2$
is treated in greater detail. The solution $u$ admits in $x$ a meromorphic  continuation
and is non-degenerate ($\det u \not= 0$). Therefore it seems
natural to call this solution fundamental even in the case of
potentials with singularities. This definition agrees with the
standard  requirement for the fundamental and Jost
solutions in the case of singularities   to be defined
by the same formula on the whole domain \cite{APP,M2}. Using (\ref{1.12})
we shall construct meromorphic in $x$  solutions $Y$ and $Z$ that satisfy (\ref{1.7}).

 Formulas (\ref{1.4}) and (\ref{1.5}) yield
\begin{equation} \label{1.14}
\a S(x)-S(x) \a^*=i \Lambda(x)j \Lambda(x)^*.
\end{equation}
(Compare relations (\ref{1.13}) and (\ref{1.14}) with the relations (\ref{0.1'}) and (\ref{0.2'}).)
 By the proof of Theorem 3.1 in \cite{SaAp} we obtain
\begin{Pn} \label{Pn1.1} Let
$\a$, $S_0$, $\g_1$, and $\g$ be an admissible set of parameter
matrices. Then we have
\begin{equation}\label{1.15}
\lim_{x \to \infty}w_{\alpha }(x,\lambda  )=
\left[\begin{array}{lr} I_{m_1}
 & 0 \\
0 &   \chi(\lambda) \end{array}\right], \quad \chi(\lambda)=
I_{m_2}-i\gamma^*\vk_R(\lambda I_n -\alpha)^{-1}\g,
\end{equation}
where
\begin{equation} \label{1.16}
 \vk_R= \lim_{x \to \infty}R(x)^{-1},
 \quad R(x)=e^{-ix \alpha } S(x)e^{ix \alpha^{*} }.
\end{equation}
\end{Pn}
Quite similar to (\ref{n1.1}) it follows from (\ref{1.13}) and (\ref{1.14}) that
\[
w_{\alpha }(x,\lambda  )^*jw_{\alpha }(x,\lambda  )=
\]
\begin{equation} \label{1.17}
=j+i(\ov{\lambda}- \lambda) \Lambda(x)^* (\ov{\lambda} I_n
-\alpha^*)^{-1}S(x)^{-1}(\lambda I_n -\alpha)^{-1}\Lambda(x).
\end{equation}
By Theorem 2.1 in \cite{SaAp} matrix function $Q(x)=e^{ix \alpha }
S(x)e^{-ix \alpha^{*} }$ is monotonic and tends to infinity, when
$x \to \infty$. Thus there is $x_0$ such that $S(x)>0$ on $(x_0,
\, \infty)$. According to (\ref{1.13}) and (\ref{1.17}) for $x \in
(x_0, \, \infty)$ the matrix function $w_{\alpha }(x,\lambda  )$
is rational in $\lambda$ and $j$-contractive in  $\BC_-$.

From Proposition \ref{Pn1.1} it follows
\begin{Pn}\label{Pn2.1}
Let $v \in $ PE be determined by the admissible matrices $\a$,
$S_0 \,$ $(\det \, S_0 \not=0)$, $\gamma_1$, and $\gamma_2$. Then
the solutions $Y(x,\lambda )$ and $Z(x,\lambda )$ of the system
$(\ref{1.1})$ defined by the boundary conditions $(\ref{1.7})$ are
given by
\begin{equation} \label{2.1}
Y(x,\lambda ) = w_{\alpha} (x,\lambda )e^{ix\lambda  j} \left[
\begin{array}{c}
I_{m_1} \\
0
\end{array}\right], \quad Z(x,\lambda )=
w_{\alpha}(x,\lambda )e^{ix\lambda  j} \left[
\begin{array}{c}
\Xi_{1}(\lambda )\\ \Xi_{2}(\lambda )
\end{array}\right],
\end{equation}
where
\begin{equation}\label{2.2}
\Xi_{1}(\lambda )=-i\g_1^*(\lambda I_n-\a^*)^{-1}S_0^{-1}
\g,\qquad \Xi_{2}(\lambda )=
\end{equation}
\[
=I_{m_2}+i\g^*(\lambda I_n-\a^*)^{-1}S_0^{-1} \g.
\]
The asymptotics of $Z(x,\lambda )$ at infinity is given by
\begin{equation}\label{2.3}
\lim_{x\to\infty}e^{-ix\lambda  j}Z(x,\lambda )=
 \left[
\begin{array}{c}
\Xi_{1}(\lambda )\\ \noalign{\vskip6pt} \chi (\lambda)
\Xi_{2}(\lambda )
\end{array}\right]\qquad (\lambda \in \BR).
\end{equation}
\end{Pn}
\begin{proof}. The proof of equalities  (\ref{2.1}) and (\ref{2.3})
is straightforward.
To prove  equalities  (\ref{2.2}) notice that
\begin{equation}\label{2.4}
\left[
\begin{array}{c}
\Xi_{1}(\lambda )\\ \Xi_{2}(\lambda )
\end{array}\right] = w_{\alpha} (0,\lambda )^{-1}\left[
\begin{array}{c} 0 \\
I_{m_2}
\end{array}\right].
\end{equation}
By (\ref{1.13}) and (\ref{1.17}) we have
\begin{equation}\label{2.5}
w_{\alpha} (x,\lambda )^{-1}=jw_{\alpha }(x, \ov{\lambda }
)^*j=I_m-i j \Lambda(x)^* (\lambda I_n
-\alpha^*)^{-1}S(x)^{-1}\Lambda(x).
\end{equation}
Formulas (\ref{2.4}) and (\ref{2.5}) yield  (\ref{2.2}).
\end{proof}
From Proposition \ref{Pn2.1} we derive
\begin{equation} \label{2.6}
Y_1(0,\lambda)=I_{m_1}+i\g_1^*S_0^{-1}(\lambda I_n-\a)^{-1}\g_1,
\quad Y_2(0,\lambda)=
\end{equation}
\[
=-i\g^*S_0^{-1}(\lambda I_n -\a)^{-1}\g_1,
\]
and
\begin{equation}
\lim_{x\to\infty}e^{-ix\lambda }Z_1(x,\lambda )=\Xi_{1}(\lambda ),
\quad \lim_{x\to\infty}e^{ix\lambda }Z_2(x,\lambda )= \chi
(\lambda) \Xi_{2}(\lambda ).\label{2.7}
\end{equation}
Notice further that, according to formula (3.16) in \cite{SaAp}, we have
\begin{equation} \label{2.7'}
\vk_R
\a- \a^* \vk_R
+i \vk_R
\g \g^*
\vk_R
=0.
\end{equation}
Thus one can easily check that
\begin{equation} \label{2.7''}
\chi(\lambda)^{-1}= I_{m_2}+i\gamma^*(\lambda I_n
-\alpha^*)^{-1}\vk_R \g.
\end{equation}
Using definitions (\ref{1.8})-(\ref{1.10}) and realizations (\ref{2.6}), (\ref{2.7}), and
(\ref{2.7''})
we obtain the next theorem
that gives explicit expressions for the reflection and transmission coefficients.

\begin{Tm} \label{Tm2.2}  Let $v\in$ PE
be determined by the admissible matrices $\a$,
$S_0 \,$ $(\det \, S_0 \not=0)$,  $\gamma_1$,  and
$\gamma$. Then the transmission and reflection coefficients are given
by the formulas
\begin{equation}
T_L(\lambda)=I_{m_1} - i\gamma_1^*S_0^{-1}(\lambda  I_n -
\theta)^{-1} \gamma_1,\qquad  \theta = \a  - i\g_1\g_1^*S_0^{-1},
\label{2.8}
\end{equation}
\begin{equation}
 R_L(\lambda)=- i\g^*S_0^{-1}(\lambda  I_n -
\theta)^{-1} \gamma_1, \label{2.9}
\end{equation}
\begin{equation}
 T_R(\lambda)= I_{m_2} - i\gamma^*S_0^{-1}(\lambda  I_n -
\theta)^{-1}
(I_n- S_0 \vk_R)\gamma , \label{2.10}
\end{equation}
\begin{equation}
 R_R(\lambda)=- i\gamma_1^*S_0^{-1}(\lambda  I_n -
\theta)^{-1}(I_n-S_0 \vk_R)\gamma- i\gamma_1^*(\lambda  I_n
-\a^*)^{-1}\vk_R\g. \label{2.11}
\end{equation}
\end{Tm}
\begin{proof}.
Compare (\ref{2.6}) with (\ref{nn1.1}) to derive (\ref{2.8}) and (\ref{2.9}) from
(\ref{nn1.2}) and (\ref{n1.3}). Notice further that in view of (\ref{0.2'}) we have
\begin{equation} \label{n3.1}
\alpha^*-i S_0^{-1} \g \g^*=S_0^{-1}(S_0 \a^*-i \g \g^*)=S_0^{-1}(\a S_0 -i \g_1 \g_1^*)=S_0^{-1}\t S_0,
\end{equation}
where $\t$ is given by (\ref{2.8}).
Therefore, taking into account Appendix one can see that
\[
\Xi_2(\lambda)^{-1}=I_{m_2}-i \g^*(\lambda I_n-(\a^*-i S_0^{-1} \g
\g^*))^{-1}S_0^{-1} \g=
\]
\begin{equation} \label{n3.2}
=
I_{m_2}-i \g^*S_0^{-1}(\lambda I_n- \t)^{-1} \g.
\end{equation}
By the identity (\ref{0.2'}) and the second relation in (\ref{2.8}) we get
\begin{equation} \label{n3.3}
 \g \g^*=i( \t S_0 - S_0 \a^*)= -i \big( (\lambda I_n- \t) S_0 - S_0 (\lambda I_n-\a^*) \big).
\end{equation}
According to (\ref{n3.2}) and (\ref{n3.3}) by direct calculation similar to the one used to derive
(\ref{n1.3}) we obtain
\begin{equation} \label{n3.4}
\Xi_1(\lambda) \Xi_2(\lambda)^{-1}= -i \g_1^*S_0^{-1}(\lambda I_n-
\t)^{-1} \g.
\end{equation}
Finally, by (\ref{1.9}), (\ref{1.10}), (\ref{2.7}), (\ref{n3.2}), and (\ref{n3.4}) we have
\begin{equation} \label{n3.5}
 R_R(\lambda)=-i \g_1^*S_0^{-1}(\lambda I_n- \t)^{-1} \g \chi(\lambda)^{-1},
\end{equation}
\begin{equation} \label{n3.6}
 T_R(\lambda)=(I_{m_2}-i \g^*S_0^{-1}(\lambda I_n- \t)^{-1} \g) \chi(\lambda)^{-1}.
\end{equation}
In view of (\ref{2.7''}) and (\ref{n3.3}), formula  (\ref{n3.6}) yields (\ref{2.10}), whereas
formula  (\ref{n3.5}) yields (\ref{2.11}).
For the case $S_0=I_m$ similar formulas have been proved in \cite{GKS6}.
\end{proof}
By (\ref{1.11}) Theorem \ref{Tm2.2} provides explicit formulas for
the scattering matrix function ${\mathcal S}(\lambda )$, too.
According to (\ref{2.8}), (\ref{2.9}), (\ref{n3.5}), and
(\ref{n3.6}) we have
\begin{equation} \label{n3.7}
{\mathcal S}(\lambda )= {\bf S}(\lambda )\left[
\begin{array}{cc} I_{m_1} & 0 \\ 0 &
\chi (\lambda )^{-1}
\end{array}
\right],
\end{equation}
where ${\bf S}(\lambda )$ is given by (\ref{nn1.7}).
%%%%%%%%%%%%%%%%%%%%%%%%%%%%%%%%%%%%%%%%%%%%%%%%%%%%%%%%%%%%%%%%%%%%%%%%%%%%%%%%%%%%%%%%%%%%%%%%%
%%%%%%%%%%%%%%%%%%%%%%%%%%%%%%%%%%%%%%%%%%%%%%%%%%%%%%%%%%%%%%%%%%%%%%%%%%%%%%%%%%%%%%%%%%%%%%%%%%%%%
\section{An inverse problem for Dirac type systems} \label{inv}
\setcounter{equation}{0}
In this section we shall consider the inverse problem, that is we shall recover system (\ref{1.1}), (\ref{1.3})
or equivalently $v \in {\mathrm{PE}}$ by $R_L$.
\begin{Tm} \label{Tm2.3}
Let $R$ be a  proper rational $m_2 \times m_1$ matrix function. Then   $R$ is the
left reflection coefficient of a system $($\ref{1.1}$)$ with $v \in {\mathrm{PE}}$
if and only if $R$
vanishes at infinity and is  contractive on $\BR$.
If $R$ satisfies these conditions, then $v$ can be uniquely
recovered from $R$ in two steps. First
recover a set of parameter matrices $\a$, $S_0$, $\g_1$, and $\g$ as described in
Theorem \ref{Tm0.1}. Next construct $v$ via formulas  $($\ref{1.3}$)$ and $($\ref{1.4}$)$.
\end{Tm}
\begin{proof}.
From (\ref{2.9}) it follows that the left reflection coefficient
is always proper rational and vanishes at infinity. Notice that
the four block scattering matrix function ${\mathcal S}(\lambda )$
is unitary on $\BR$. Indeed, by (\ref{n3.7}) ${\mathcal S}(\lambda
)$ is the product of two matrix functions. The first factor ${\bf
S}(\lambda )$ is unitary by Corollary \ref{Cy1.1}, and the second
one is unitary because it is $j$-unitary (see formula
(\ref{1.15})) and diagonal. As  ${\mathcal S}(\lambda )$ is
unitary, so its block $R_L$ is contractive, and the necessary
conditions are proved.

 Now let $R$ be a proper rational $m_2 \times
m_1$ matrix function which vanishes at infinity and is
contractive on $\BR$. Then according to the proof of Theorem \ref{Tm0.1}
the identity (\ref{1.5}) holds for the recovered by $R$ parameter matrices $\a$, $S_0$, $\g_1$, and $\g$.
Moreover, $S_0$ is invertible and in view of  (\ref{2.16}) the first relation in (\ref{1.6}) is true.
Recall (see also Appendix)  that
if a pair $A$, $B$ is full range,  then
the pair $A-BK$, $B$ is also full range. By the minimality of the
realization (\ref{2.15}) the pairs $A$, $B$ and $A^*$, $C^*$ are full range. Thus by (\ref{0.4}) the pair
$\a$, $\g_1$  is full range,  i.e., we get the second relation in (\ref{1.6}). Finally, notice
that the pair $S_0 A^*S_0^{-1}$, $S_0C^*$ is full range simultaneously with the pair $A^*$, $C^*$.
By (\ref{0.4}), (\ref{nn1.6}), and (\ref{n3.1}) it means that  the pair $\a + i \g \g^*S_0^{-1}$, $\g$ is
full range, and so the pair $\a $, $\g$ is
full range too. Therefore the third relation in (\ref{1.6}) is true as well. According to
(\ref{1.5}) and (\ref{1.6}), the set $\a$, $S_0$, $\g_1$, and $\g$
is admissible.  Hence we can apply Theorem \ref{Tm2.2}. As the realization
(\ref{2.15}) can be rewritten in the form
\[
R(\lambda)=- i\g^*(\lambda  I_n - \theta)^{-1} \gamma_1,
\]
where $\theta =A= \a - i \g_1 \g_1^*S_0^{-1}$, one can see that $R$ is the left
reflection coefficient $R_L$ of the system (\ref{1.1}) with the
potential $v$ determined by the admissible set $\a$, $S_0$, $\g_1$, and $\g$.

To prove the uniqueness of  the recovered $v$ suppose that there is another system with
the potential $\wt v$ determined by the admissible set $\wt \a$, $\wt S_0$, $\wt \g_1$,
and $\wt \g$, and with the same left reflection coefficient $R$. Repeat the arguments in the proof
of the uniqueness in Theorem \ref{Tm0.1} to get
\begin{equation}\label{nn4.1}
\wt S_0=sS_0s^*, \quad \wt \Lambda_0=s \Lambda_0, \quad \wt \a = s
\a s^{-1}.
\end{equation}
Taking into account (\ref{1.4}) and (\ref{nn4.1}) we derive
\begin{equation}\label{nn4.2}
\wt \Lambda(x)=s \Lambda(x), \quad \wt S(x)=sS(x)s^*,
\end{equation}
where $\wt S$ and $\wt \Lambda$ are the matrix functions defined
via (\ref{1.4}) for the second admissible set. Then formulas
(\ref{1.3}), (\ref{nn4.1}), and the second equality in
(\ref{nn4.2}) yield $v= \wt v$.
\end{proof}

\section{Appendix on mathematical system theory} \label{ap}
\setcounter{equation}{0}
The material from  mathematical system theory of rational matrix
functions, that is used in this paper, has its roots in the Kalman
theory  \cite{KFA}, and can be found in various books
(see, for instance,  \cite{BGK,CF,LR}). Our presentation closely follows
the presentation in \cite{GKS6}. The rational
matrix functions appearing in this paper are {\em proper}, that
is, analytic at infinity.  Such an $m_2 \times m_1$ matrix function $F$
can be represented in the form
\begin{equation}
\label{app.1} F(\lambda)=D+C(\lambda I_n-A)^{-1}B,
\end{equation}
where $A$ is a square matrix of some order $n$,
the matrices $B$ and $C$ are of sizes $n \times m_1$
and $m_2 \times n$, respectively, and $D=F(\infty)$.
The representation (\ref{app.1}) is called
a {\it realization} or a {\it transfer matrix representation} of
$F$, and the number $\mathrm{ord}(A)$ (order of the matrix $A$) is called the {\it state
space dimension} of the realization.

 The realization
(\ref{app.1}) is said to be {\it minimal} if its state space
dimension $n$ is minimal among all possible realizations of $F$.
This minimal $n$  is
called the {\it McMillan degree} of $F$.
The realization (\ref{app.1}) of $F$ is minimal if and only if
\begin{equation}
\label{app.2'} {\mathrm{span}}\bigcup_{k=0}^{n-1}\im
A^kB=\BC^n,\quad \bigcup_{k=0}^{n-1}\im  (A^*)^kC^*=\BC^n,
\quad n= \mathrm{ord}(A).
\end{equation}
If for a pair of matrices $A$, $B$ the first equality in
(\ref{app.2'}) holds, then the pair $A$, $B$ is called {\it
controllable} or a {\it full range}. If the second equality
in (\ref{app.2'}) is fulfilled, then $C$, $A$ is said to be
{\it observable}. If a pair $A$, $B$
 is full range, and $K$ is an $m_1 \times n$ matrix,  then
the pair $A-BK$, $B$ is also full range.

Minimal realizations are unique up to a basis transformation, that
is, if (\ref{app.1}) is a minimal realization of $F$ and if
$F(\lambda)=D+\widetilde C(\lambda I_n -\widetilde
A)^{-1}\widetilde B$ is a second minimal realization of $F$, then
there exists an invertible matrix $s$ such that
\begin{equation}
\label{app.3} \widetilde A=sAs^{-1},\quad \widetilde B=sB,\quad
\widetilde C=Cs^{-1}.
\end{equation}
In this case, (\ref{app.3})  is called a {\it similarity transformation}.

Finally, if $F$ is a square matrix and $D=I_m$, then $F^{-1}$ admits representation
\begin{equation} \label{app4}
F(\lambda)^{-1}=I_m-C(\lambda I_n-A^\times)^{-1}B,\quad
A^\times=A-BC.
\end{equation}

{\bf Acknowledgment.} AS is grateful to the German Academic
Exchange Service (DAAD) for the support of this work through grant
A/04/15892 and to the Leipzig University, where this work was
fulfilled.

\newpage
Bernd Fritzshe \\ Fakult\"at f\"ur Mathematik und Informatik,
Mathematisches Institut, Universit\"at Leipzig, Augustusplatz
10/11, D-04109 Leipzig, Germany;
\\
Bernd Kirstein \\ Fakult\"at f\"ur Mathematik und Informatik,
Mathematisches Institut, Universit\"at Leipzig, Augustusplatz
10/11, D-04109 Leipzig, Germany; \\ Alexander Sakhnovich \\ Branch
of Hydroacoustics,  Marine Institute of Hydrophysics, National
Academy of Sciences, Preobrazhenskaya 3,  Odessa, Ukraine. \\


\begin{thebibliography}{FKS}
\bibitem{AKNS}
M.J. Ablowitz, D.J. Kaup, A.C. Newell, H. Segur,  The inverse
scattering transform - Fourier analysis for nonlinear problems,
Stud. Appl. Math. 53 (1974) 249--315.

\bibitem{AS}
M.J. Ablowitz, H. Segur,   Solitons and the inverse
scattering transform, in: SIAM Stud. Appl.
Math.,  Vol. 4, Philadelphia, 1981.

\bibitem{AG0}
D. Alpay, I. Gohberg, Unitary rational matrix functions, in: Oper. Theory Adv. Appl., Vol.33, 175-222,
Birkh\"auser, Basel, 1988.

\bibitem{AG}
D. Alpay, I. Gohberg, Inverse spectral problem for
differential operators with rational scattering matrix functions,
J. Diff. Eqs  118 (1995) 1--19.

\bibitem{AGKS}
D. Alpay, I. Gohberg, M.A. Kaashoek, A.L.  Sakhnovich,  Direct
and inverse scattering problem for canonical systems with a
strictly pseudo-exponential potential, Math. Nachr. 215 (2000)
5--31.

\bibitem{AKM}
T. Aktosun, M. Klaus, C. van der Mee, Direct and inverse
scattering for selfadjoint Hamiltonian systems on the line,
Integral Equations Operator Theory
38 (2000) 129--171.




\bibitem{APP}
V.A. Arkadiev, A.K. Pogrebkov, M.K. Polivanov, Singular
solutions of the KdV equation and the method of the inverse
problem (in Russian),    Nauchn. Sem. Leningrad.
Otdel. Mat. Inst. Steklov. (LOMI)  133 (1984)  17--37.

\bibitem{Ar}
D.Z. Arov,
Realization of matrix-valued functions according to Darlington (Russian),
Izv. Akad. Nauk SSSR Ser. Mat. 37 (1973) 1299--1331.

\bibitem{Ar6}
D.Z.Arov,  Gamma-generating matrices, $j$-inner matrix functions and related extrapolation problems I,
J. Soviet Math. 52 (1990) 3487--3491.

\bibitem{AD}
D.Z. Arov, H. Dym, $J$-inner matrix functions,
interpolation and inverse problems for canonical systems V,
Integral Equations Operator Theory  43 (2002) 68--129.

\bibitem{AFK0}
D.Z. Arov, B. Fritzsche, B. Kirstein,
Completion problems for $j\sb {pq}$-inner functions I, Integral Equations Operator Theory 16 (1993) 155--185.

\bibitem{AFK}
D.Z. Arov, B. Fritzsche, B. Kirstein,
Completion problems for $j\sb {pq}$-inner functions II, Integral Equations Operator Theory 16 (1993) 453--495.

\bibitem{B1}
J.A. Ball, E.A. Jonckheere, The four-block Adamjan-Arov-Krein problem,
J. Math. Anal. Appl. 170 (1992) 322-342.


\bibitem{B2}
J.A. Ball, Linear systems, operator model theory and scattering: multivariable generalizations,
Fields Inst. Commun. 25 (2000) 151--178.

\bibitem{BKo}
S. Barran, M. Kovalyov, A note on slowly decaying solutions
of the defocusing nonlinear Schr\"odinger equation,   J. Phys.A
32 (1999) 6121--6125.

\bibitem{BGK}
H. Bart,  I. Gohberg,  M.A. Kaashoek,  Minimal  factorization of
matrix  and  operator  functions,  Birkh\"auser, Basel, 1979.


\bibitem{BM}
L.A. Bordag, V.B.  Matveev, Selfsimilar solutions of the
Korteweg-de Vries equation and potentials with a trivial
$S$-matrix,   Theoret. and Math. Phys. 34 (1978) 272--275.


\bibitem{BC}
R. Beals, R.R.  Coifman, Scattering and inverse scattering
for first order systems,  Comm. Pure Appl. Math. 37 (1984)
39--90.

\bibitem{BDT}
R. Beals, P. Deift, C. Tomei,  Direct and inverse
scattering on the line, in: Mathematical Surveys and Monographs, Vol.
28, AMS, Providence,  1984.


\bibitem{CS} K. Chadan, P.C.  Sabatier,  Inverse problems in
quantum scattering theory, Springer, NY, 1989.

\bibitem{CG} S. Clark, F.  Gesztesy,
Weyl-Titchmarsh $M$-function asymptotics, local uniqueness
results, trace formulas, and Borg-type theorems for Dirac
operators, Trans. Amer. Math. Soc. 354 (2002) 3475--3534.

\bibitem{CF}
M.J. Corless,  A.E. Frazho,  Linear Systems and Control - An
Operator Perspective, Marcel Dekker, 2003.

\bibitem{DeD}
P. Dewilde, H. Dym,  Lossless chain scattering matrices
and optimum linear prediction: the vector case, Internat. J. Circuit
Theory Appl. 9 (1981) 135--175.


\bibitem{DFK}
V.K. Dubovoj, B. Fritzsche, B. Kirstein,
Matricial version of the classical Schur problem, in:
Teubner-Texte zur Mathematik [Teubner Texts in Mathematics], Vol. 129,
B. G. Teubner Verlagsgesellschaft mbH, Stuttgart, 1992.

\bibitem{DI}
H. Dym, A. Iacob, Positive definite extensions, canonical
equations and inverse problems,  in: Oper. Theory Adv. Appl., Vol.  12,
141-240,
 Birkh\"auser, Basel-Boston, 1984.

\bibitem{FT}
L.D. Faddeev, L.A. Takhtajan,   Hamiltonian methods in
the theory of solitons, Springer, NY, 1986.

\bibitem{FK1}
B. Fritzsche, B. Kirstein, M. Mosch,
On block completion problems for Arov-normalized
$j\sb {qq}$-$J\sb q$-elementary factors,
Lin. Alg. Appl. 346 (2002) 273--291.

\bibitem{FK2}
B. Fritzsche, B. Kirstein, K. M\"uller,
An analysis of the block structure of $j\sb {qq}$-inner functions,
 in: Oper. Theory Adv. Appl., Vol. 106, 157--185, Birkh\"auser, Basel, 1998.

\bibitem{GSS}
F. Gesztesy, W.  Schweiger, B. Simon, Commutation methods
applied to the mKdV-equation, Trans. Amer. Math. Soc.
324 (1991) 465--525.



\bibitem{GKS1}
I. Gohberg I,  M.A. Kaashoek, A.L.  Sakhnovich,  Canonical
systems with rational spectral densities: explicit formulas and
applications, Math. Nachr. 194 (1998) 93--125.

\bibitem{GKS4}
I. Gohberg, M.A.  Kaashoek,  A.L.  Sakhnovich, Canonical
systems on the  line with rational spectral densities: explicit
formulas,  in: Oper. Theory Adv. Appl., Vol. 117, 127-139, Birkh\"auser,
Basel, 2000.

\bibitem{GKS6}
I. Gohberg, M.A.  Kaashoek,  A.L.  Sakhnovich, Scattering
problems for a canonical system with a pseudo-exponential
potential,  Asymptotic Analysis 29 (2002) 1--38.

\bibitem{GR}
I. Gohberg, S. Rubinstein,
Proper contractions and their unitary minimal completions,
in: Oper. Theory Adv. Appl., Vol. 33, 223--247, Birkh\"auser, Basel, 1988.

\bibitem{KFA}
R.E. Kalman, P.L. Falb,  M.A. Arbib,  Topics in
mathematical system theory, McGraw-Hill, New York, 1969.

\bibitem{K}
M.G. Krein, Topics in differential and integral equations
and operator theory,  Oper. Theory Adv. Appl., Vol. 7,
Birkh\"auser, Basel-Boston, 1983.

\bibitem{Kur}
P. Kurasov, Singular and supersingular perturbations: Hilbert space
methods. Spectral theory of Schrödinger operators,  Contemp. Math.
 340 (2004) 185-216.

\bibitem{LR}
P. Lancaster, L. Rodman,
 Algebraic Riccati equations, Clarendon Press, Oxford,
1995.

\bibitem{LS}
B.M. Levitan, I.S.  Sargsjan,  Sturm-Liouville and Dirac
operators, in: Mathematics and its Applications, Vol. 69,
Kluwer, Dordrecht, 1991.


\bibitem{L}
M.S. Liv\v{s}ic,  On a class of linear operators in Hilbert
space, in:   Amer. Math. Soc. Transl. (2), Vol. 13, 85-103, 1960.

\bibitem{M}
V.A. Marchenko, Nonlinear equations and operator algebras,
 Reidel, Dordrecht, 1988.

\bibitem{M1}
V.B. Matveev,  Asymptotics of the multipositon-soliton
$\tau$ function of the Korteweg-de Vries equation and the supertransparency,
J. Math. Phys. 35 (1994) 2955--2970.

\bibitem{M2}
V.B. Matveev,  Positons: slowly decreasing analogues of solitons,
Theoret. and Math. Phys. 131 (2002) 483--497.



\bibitem{SaA1}
A.L. Sakhnovich, Exact solutions of nonlinear equations and the
method of operator identities,  Lin.Alg.Appl. 182 (1993)
109--126.


\bibitem{SaA3}
A.L. Sakhnovich, Dirac type and canonical  systems: spectral
and Weyl-Titchmarsh fuctions, direct and inverse problems,
Inverse Problems  18 (2002) 331--348.

\bibitem{SaAp}
A.L. Sakhnovich, Dirac type system on the axis: explicit
formulas for matrix potentials with singularities and
soliton-positon interactions, Inverse Problems  19 (2003),
845-854.

\bibitem{SaL1}
L.A. Sakhnovich, On  the  factorization  of  the  transfer
matrix function,  Sov. Math. Dokl. 17 (1976) 203--207.


\bibitem{SaL2}
L.A. Sakhnovich, Factorisation  problems  and  operator
identities,   Russian
 Math. Surv. 41 (1986) 1--64.


\bibitem{SaL3}
L.A. Sakhnovich, Spectral theory of canonical differential
systems. Method of operator identities, in: Oper. Theory Adv. Appl., Vol.107,
Birkh\"auser, Basel-Boston, 1999.

\bibitem{Sch} C. Schiebold, Solitons of the sine-Gordon
equation coming in clasters,  Revista Mat. Compl. 15 (2002)
262--325.

\bibitem{V}
V. Vinnikov,
Commuting operators and function theory on a Riemann surface,
 in:
Math. Sci. Res. Inst. Publ., Vol. 33, 445--476,
Cambridge Univ. Press, Cambridge, 1998.


\bibitem{YL} A.E. Yagle, B.C.  Levy,
The Schur algorithm and its applications,   Acta Appl. Math.
 3 (1985) 255--284.

\bibitem{ZS}
V.E. Zaharov, A.B. Shabat, On soliton interaction in stable
media,  JETP 64 (1973) 1627--1639.

\bibitem{Zh}
X. Zhou, Inverse scattering transform for systems with rational
spectral dependence, J. Diff. Eqs.  115 (1995)  277--303.


% \bibitem{label}
% Text of bibliographic item

% notes:
% \bibitem{label} \note

% subbibitems:
% \begin{subbibitems}{label}
% \bibitem{label1}
% \bibitem{label2}
% If there is a note, it should come last:
% \bibitem{label3} \note
% \end{subbibitems}



\end{thebibliography}
\end{document}